\documentclass[12pt,oneside,a4paper]{article}

\title{Fixed point sets of equivariant fiber-preserving maps}
\author{ \normalsize Rafael Souza \\ 
        \small Universidade Estadual de Mato Grosso do Sul\\
				\small  moreira@uems.br \\
				\and 
				\normalsize Peter Wong \\
				\small Bates College\\
				\small  pwong@bates.edu
				}
\date{\today}

\usepackage{amsmath}
\usepackage{amsfonts}
\usepackage{amssymb}
\usepackage{amsthm}
\usepackage[T1]{fontenc}
\usepackage[all]{xy}
\usepackage{graphicx}
\usepackage[english]{babel}

\theoremstyle{definition}
\newtheorem{teo}{Theorem}[section]
\newtheorem{lema}[teo]{Lemma}
\newtheorem{prop}[teo]{Proposition}
\newtheorem{cor}[teo]{Corollary}

\newtheorem{exem}[teo]{Example}
\newtheorem{rem}[teo]{Remark}

\begin{document}

\maketitle

{\bf keywords:} fixed points, equivariant maps

{\bf AMS subj. class. [2010]:} Primary: 55M20; Secondary: 57S99

\begin{abstract}

Given a selfmap $f:X\to X$ on a compact connected polyhedron $X$, H. Schirmer gave necessary and sufficient conditions for a nonempty closed subset $A$ to be the fixed point set of a map in the homotopy class of $f$. R. Brown and C. Soderlund extended Schirmer's result to the category of fiber bundles and fiber-preserving maps. The objective of this paper is to prove an equivariant analogue of Brown-Soderlund theorem result in the category of $G$-spaces and $G$-maps where $G$ is a finite group.
 
\end{abstract}

\section{Introduction and statement of results}

A well-known and important question in classical topology is the {\it fixed point property}. Recall that a topological space $X$ is said to have the fixed point property if every (continuous) map $f:X\to X$ must have a fixed point $x_0\in X$ such that $f(x_0)=x_0$. A related question is the so-called {\it complete invariance property for deformation} (CIPD). We say that $X$ has the CIPD if for any nonempty closed subset $A\subset X$, there exists a selfmap $f:X\to X$ homotopic to the identity $1_X$ such that $A={\rm Fix} (f)=\{x\in X \mid f(x)=x\}$. In \cite{sc1}, H. Schirmer generalized the concept of CIPD and gave necessary and sufficient conditions for a nonempty closed subset $A$ to be the fixed point set of a map $g$ in the homotopy class of a given selfmap $f$. That is, given a map $f:X\to X$, Schirmer determined when a closed nonempty subset $A$ can be realized as $A={\rm Fix} (g)$ for some $g$ homotopic to $f$. Upon relaxing the conditions given by Schirmer, C. Soderlund together with R. Brown \cite{bs} generalized Schirmer's result to fiber-preserving maps of fiber bundles.

Suppose that $X$ is a compact connected polyhedron without local cutpoints and $A$ is a closed subset imbedded inside a subpolyhedron $K$ that can be {\it by-passed in $X$}, that is, every path $C$ in $X$ with $C(0), C(1) \in X-K$, is homotopic to a path $C'$ in $X-K$ relative to the endpoints. H. Schirmer \cite{sc1} introduced the following two conditions to realize $A$ as the fixed point set of a selfmap in the homotopy class $[f]$.

\begin{enumerate}
\item[(C1)] if there exists a homotopy $H_A:A\times [0,1]\to X$ from $f|_A$ to the inclusion $i:A\hookrightarrow X$;
\item[(C2)] if for every essential fixed point class $\mathbb F$ of $f$, there exists a path $\alpha :[0,1]\to X$ with $\alpha(0)\in \mathbb F, \alpha(1)\in A$ and $\{\alpha(t)\}\sim \{f\circ \alpha(t)\} \ast \{H_A(\alpha(1),t)\}$ relative to the endpoints.
\end{enumerate}

Soderlund \cite[Theorem 3.5]{s} showed, by relaxing the assumption on $A$ given by Schirmer, the following result.

\begin{teo}\label{soderlund} Let $X$ be a compact, connected polyhedron with no local cut points and $A$ be a closed locally contractible subspace of $X$ such that $X-A$ is not a $2$-manifold and $A$ can be by-passed in $X$. Then $A={\rm Fix} (g)$ for some $g\sim f$ if and only if (C1) and (C2) are satisfied.
\end{teo}

Subsequently, R. Brown and C. Soderlund \cite{bs} introduced analogous conditions in the fiber-preserving setting. Let $\mathfrak F=(E,p,B;Y)$ be a (locally trivial) fiber bundle and $f:E\to E$ a fiber preserving map.

\begin{enumerate}\label{fiber-conditions}
\item[(${\rm C1}_{\mathfrak F}$)] if there exists a fiber preserving homotopy $H_A:A\times [0,1]\to E$ from $f|_A$ to the inclusion $i:A\hookrightarrow E$;
\item[(${\rm C2}_{\mathfrak F}$)] if for every essential fixed point class $\mathbb F$ of $f$, there exists a path $\alpha :[0,1]\to E$ with $\alpha(0)\in \mathbb F, \alpha(1)\in A$ and $\{\alpha(t)\}\sim \{f\circ \alpha(t)\} \ast \{H_A(\alpha(1),t)\}$ relative to the endpoints.
\end{enumerate}

Following the terminology of \cite{bs}, we call $(X,A)$ a {\it suitable pair} if $X$ is a finite polyhedron with no local cut points and $A$ is a closed locally contractible subspace of $X$ such that $X-A$ is not a $2$-manifold and $A$ can be by-passed in $X$.

In \cite {bs}, it was shown that conditions $({\rm C1}_{\mathfrak F})$ and $({\rm C2}_{\mathfrak F})$ are also sufficient. The following is their main result.

\begin{teo}\label{brown-soderlund}
Let $\mathfrak F=(E,p,B;Y)$ be a fiber bundle where $E,B$ and $Y$ are connected finite polyhedra, $f:E\to E$ a fiber preserving map and $A$ a closed locally contractible sub-bundle of $E$ such that each component $p(A)_j$ of $p(A)$ is contractible and $(B,p(A)), (Y,Y_j)$ for all sub-bundle fibers $Y_j$ of $A$, are suitable pairs. Suppose $({\rm C1}_{\mathfrak F})$ and $({\rm C2}_{\mathfrak F})$ are satisfied and $A$ intersects every essential fixed point class of $f_{b_j}:p^{-1}(b_j)\to p^{-1}(b_j)$ for at least one $b_j$ in each component $p(A)_j$. If $Z$ is a closed bundle subset of $A$ that intersects every component of $A$, then there exists a map $g:E\to E$ that is fiber preserving and fiberwise homotopic to $f$ ($g\sim_{\mathfrak F} f$) such that ${\rm Fix} (g)=Z$.
\end{teo}

In particular, when $Z=A$, this theorem shows that $({\rm C1}_{\mathfrak F})$ and $({\rm C2}_{\mathfrak F})$ are necessary and sufficient for $A={\rm Fix}( g)$ for some $g\sim_{\mathfrak F} f$.

Many applications involve symmetries in the presence of a group action. As a result, equivariant topology has been proven to be useful in the study of nonlinear problems. In the equivariant setting, we are concerned with a group $G$ acting on a space $X$ together with a $G$-map $f:X\rightarrow X$ which respects the group action, that is, for all $\alpha\in G$, $f(\alpha x )=\alpha f(x)$ for all $x\in X$. In this case, the fixed point set ${\rm Fix}(f)$ is {\it a priori} a $G$-invariant subset of $X$.

In \cite{sc1}, Schirmer observed that for a given selfmap $f:\mathbb{S}^{n}\rightarrow \mathbb{S}^{n}$ of an $n$-sphere, $n\geq 2$, {\it any} closed nonempty proper subset $A$ of $\mathbb{S}^{n}$ can be realized as the fixed point set of a map $g\in[f]$ with ${\rm Fix}(g)=A$. However, such phenomenon does not hold if we impose a group action as we show in the following example, which gives the underlying motivation for this paper.

\begin{exem}\label{motivation}
Let $G=\mathbb{Z}_{2}$, $X=\mathbb{S}^{2}$ and the action  is given by $\xi(x,y,z)\mapsto(-x,-y,z)$. If $A=\{(x,y,0)\in \mathbb{S}^{2}\}$ then $A$ is $\mathbb{Z}_{2}$-invariant, but there is no $\mathbb{Z}_{2}$-map $h:\mathbb{S}^{2}\rightarrow\mathbb{S}^{2}$ that is  $\mathbb{Z}_{2}$-homotopic to the identity map $Id:\mathbb{S}^{2}\rightarrow \mathbb{S}^{2}$ such that ${\rm Fix}(h)=A$.

In fact, suppose  there is a $\mathbb{Z}_{2}$-homotopy $H$ from $Id$ to $h:\mathbb{S}^{2}\rightarrow \mathbb{S}^{2}$ such that ${\rm Fix}(h)=A$. Then, $h$ preserves $X^{G}=\{N,S\}$, where $N=(0,0,1)$ and $S=(0,0,-1)$. Hence, $h(N)=S$ and $h(S)=N$ and the path $p:I\rightarrow X^{G}$ defined by $p(t)=H^{G}(N,t)$ is such that $p(0)=Id(N)=N$ and $p(1)=h(N)=S$. But, this is impossible.

In this situation, the location of $A$ in $X$ is more important than its topology, because if we replace $A$ by $A'=\{(x,0,z)\in \mathbb{S}^{2}\}$ then:

\begin{displaymath}
H(t, \cos\theta\sin\psi, \sin\theta\sin\psi ,  \cos\psi) = 
(\cos(\theta+t\epsilon \sin\theta)\sin\psi, \sin(\theta+t\epsilon \sin\theta)\sin\psi ,  \cos\psi),
\end{displaymath}

\noindent is a $\mathbb{Z}_{2}$-homotopy (with polar coordinates) between the identity and the $\mathbb{Z}_{2}$-map $h$ such that ${\rm Fix}(h)=A'$.
\end{exem}

\begin{exem}
It is easy to see, by modifying the last example, that the equivariant analogue of Schirmer's result does not hold in general. Let $G=\mathbb{Z}_{2}$, $X=\mathbb{S}^{2}\times \mathbb{S}^2$ and the action  is given by $\xi((x,y,z), (x',y',z'))\mapsto((-x,-y,z),(1,0,0))$. The set $A=\{((x,y,0),(1,0,0))\in X\}$ is $\mathbb Z_2$-invariant and $X^G=\{N,S\}\times \{(1,0,0)\}$ consists of two points. The same argument as in Example \ref{motivation} shows that $A$ cannot be the fixed point set of any map $\mathbb Z_2$-homotopic to the identity map while $(X,A)$ satisfies the conditions of Schirmer's result for $A$ can be by-passed since $A$ has codimension $3$ in $X$.
\end{exem}

The main objective of this paper is to give an equivariant analogue of Schirmer's result and of Brown-Soderlund's result. This paper is organized as follows. In the first section, we briefly recall the non-equivariant results of \cite{sc1} and \cite{bs} and review some basic background on $G$-maps and $G$-spaces where $G$ denotes a compact Lie group. Then we review the necessary equivariant Nielsen fixed point theory from \cite{w1}. In section 2, we prove our first main result, an equivariant analogue of \cite{sc1}:

\begin{teo}\label{Cor,sod} Let $G$ be a compact Lie group, $X$ be a compact and smooth $G$-manifold and $A$ be a nonempty, closed, locally contractible $G$-subset of $X$ such that for each finite $WK$ we assume that $\dim(X^{K})\geq 3$, $\dim(X^{K})-\dim(X^{K}-X_{K})\geq2$ and $A^{K}$ is by-passed in $X^{K}$, for all $(K)\in {\rm Iso}(X)$. Suppose that the following conditions holds for a $G$-map $f:X\rightarrow X$:

\begin{itemize}
\item[($C_G1$)]  there exists a $G$-homotopy $H_{A}:A\times I\rightarrow X$ from $f|_{A}$ to the inclusion $i:A\hookrightarrow X$;
\item[($C_G2$)]  for each finite $WK$, for every $WK$-essential fixed point class $F$ of $f^{K}:X^{K}\rightarrow X^{K}$ there exists a path $\alpha: I\rightarrow X^{K}$ with $\alpha(0)\in F$, $\alpha(1)\in A^{K}$, and $\{\alpha(t)\}\sim \{f^{K}\circ\alpha(t)\}\ast\{H_{A}^{K}(\alpha(1),t)\}.$
\end{itemize}

Then for every closed $G$-subset $\Phi$ of $A$ that has nonempty intersection with every component of $A$ there exists a $G$-map $h:X\rightarrow X$, $G$-homotopic to $f$ with ${\rm Fix}(h)=\Phi$.
\end{teo}

\noindent In the last section, we apply Theorem \ref{Cor,sod} to prove an equivariant analogue of \cite{bs} when $G$ is finite:

\begin{teo}\label{Teo;bs}
 Let $G$ be a finite group, $\mathfrak{F}=(X,p,B,Y)$ be a $G$-fiber bundle where $X$, $B$ and $Y$ are compact and smooth $G$-manifolds, $\dim(B^{K})\geq 3$, $\dim(B^{K})-\dim(B^{K}-B_{K})\geq2$, for all $(K)\in {\rm Iso}(B)$, $\dim(Y^{K})\geq 3$, $\dim(Y^{K})-\dim(Y^{K}-Y_{K})\geq2$, for all $(K)\in {\rm Iso}(Y)$.

Let $A$ be a nonempty, closed, locally contractible $G$-subset of $X$ such that $(X,A)$ is $G$-fiber bundle pair with respect to the fiber bundle $\mathfrak{F}$, $p(A)$ be a closed $G$-subset of $B$ such that each component $p(A)_{j}$ of $p(A)$ is equivariantly contractible and $p^{K}(A^{K})$ is by-passed in $B^{K}$, for all $(K)\in {\rm Iso}(B)$. Let $Y_{j}$ be a subbundle fiber of $A$ such that $Y_{j}$ is a closed and locally contractible $G$-subset of $Y$ and $Y^{K}_{j}$ is by-passed in $Y^{K}$, for all $(K)\in {\rm Iso}(Y)$, and $f:X\rightarrow X$ be a $G$-fiber-preserving map such that $A^{K}$ intersects every essential $WK$-fixed point class of $f^{K}_{b_{j}}:WK(p^{K})^{-1}(\{b_{j}\})\rightarrow WK(p^{K})^{-1}(\{b_{j}\})$ for at least one $b_{j}$ in each component $p^{K}(A^{K})_{j}$, for all $(K)\in {\rm Iso}(X)$. Suppose that the following conditions hold for $f$ and $A$:

\begin{itemize}
\item[$(C_G1)_{\mathfrak{F}}$]  there exists a $G$-fiberwise-homotopy $H_{A}:A\times I\rightarrow X$ from $f|_{A}$ to the inclusion $i:A\hookrightarrow X$;
\item[$(C_G2)_{\mathfrak{F}}$]  for every $WK$-essential fixed point class $F$ of $f^{K}:X^{K}\rightarrow X^{K}$ there exists a path $\alpha: I\rightarrow X^{K}$ with $\alpha(0)\in F$, $\alpha(1)\in A^{K}$, and $\{\alpha(t)\}\sim \{f^{K}\circ\alpha(t)\}\ast\{H_{A}^{K}(\alpha(1),t)\}.$
\end{itemize}

Then for every nonempty closed $G$-bundle subset $\Phi$ of $A$ that intersects every component of $A$ there exists a $G$-fiber-preserving map $h$, $G$-fiberwise homotopic to $f$ with ${\rm Fix}(h)=\Phi$.  
\end{teo}

In order to establish the notations, let $G$ be a topological group and $X$ be a (left) $G$-space. Given a subgroup $K$ of $G$ we denote by $NK$ the normalizer of $K$ in $G$, $WK =\frac{NK}{K}$ is the Weyl group of $K$ in $G$. The orbit type of $K$ is the conjugacy class of $K$ in $G$ denoted by $(K)$. If $(K_{1})$ is subconjugate to $(K_{2})$, we write $(K_{1})\leq(K_{2})$.

If $x\in X$, then $G_{x}=\{g\in G; \  gx=x\}$ denotes the isotropy subgroup of $x\in X$, and $(G_{x})$ is called an isotropy type of $X$. We denote by ${\rm Iso}(X)$ the set of isotropy types of $X$. Moreover, $X^{K}=\{ x\in X; K\leq G_{x}\}$, $X^{(K)} =\{ x\in X; (K)\leq (G_{x})\}$, $X_{K} =\{ x\in X; G_{x}=K\}$ and $X_{(K)} =\{ x\in X; G_{x}\subset(K)\}$.

If ${\rm Iso}(X)$ is finite (in particular when $G$ is finite), we can choose an admissible ordering on ${\rm Iso}(X)$ such that $(K_{i})\leq(K_{j})$ implies $i\leq j$. Then we have a filtration of $G$-subspaces $X_{1}\subset\dots\subset X_{n}=X$ where $X_{i}=\{x\in X;\ (G_{x})=(H_{j})$ for some $ j\leq i\}$

If $f:X\to X$ is a $G$-map, then $f^K=f|_{X^K}:X^K\to X^K$ is a $WK$-map. Let $\mathcal F=\{(K)\in {\rm Iso}(X)\mid |WK|<\infty\}$ and $(K)\in \mathcal F$. If $x,y\in {\rm Fix} (f^K)$ then $x\sim_K y$ if either $y=\alpha x$ for some $\alpha \in WK$ or $\exists \ \sigma:[0,1]\to X^K$ such that $\sigma \sim f^K\circ \sigma$ relative to endpoints. Then $\sim_K$ is an equivalence relation on ${\rm Fix} (f^K)$ and the equivalence classes are called the $WK-${\it fixed point classes} of $f^K$. Evidently, a $WK$-fixed point class $\mathcal W$ is a disjoint union of a finite number of ordinary fixed point classes $W_1,...,W_r$ of $f^K$ and thus the fixed point index $ind(\mathcal W)$ is defined as $ind(\mathcal W)=\sum_i ind(W_i)$. A $WK$-fpc (fixed point class) $\mathcal W$ is {\it essential} if $ind(\mathcal W)\ne 0$. For further information on equivariant Nielsen fixed point theory, see \cite{w1}. Throughout, by a smooth $G$-manifold $X$, we assume that the fixed point set $X^H$ is a smooth connected submanifold for each isotropy subgroup $H\le G$.

\section[Proof of Theorem 1.5]{Proof of Theorem \ref{Cor,sod} - An equivariant analogue of a result of Soderlund-Schirmer}\label{sec2}

If $X$ is a smooth $G$-manifold and $A$ is a closed smooth $G$-submanifold of $X$, $G$ being a finite group, then there exists a smooth equivariant triangulation $f_{1}:(K,K_{0})\rightarrow (X,A)$ as proved in \cite{si}.
If $B$ is another closed smooth $G$-submanifold of $X$ then there is a smooth equivariant triangulation $f_{2}:(L,L_{0})\rightarrow (X,B)$ and $G$-subdivisions $K'$ of $K$ and $L'$ of $L$ such that $f'^{-1}_{1}\circ f'_{2}:|L'|\rightarrow |K'|$ is a simplicial $G$-homeomorphism, where $f'_{1}$ and $f'_{2}$ are smooth $G$-triangulations (see \cite{si}).

By Corollary 3.3.5 of \cite{sp} and $G$ being finite, we can find unique $G$-subcomplexes $L'_{0}$ of $L$ and $K'_{0}$ of $K$ such that $L'_{0}$ is a refinement of $L_{0}$ and $K'_{0}$ is a refinement of $K_{0}$. Then, $f'^{-1}_{1}\circ f'_{2}(L'_{0})=K_{1}$ is a $G$-subcomplex of $K'$ and a $G$-triangulation of $B$. In fact, $f'_{1}(K_{1})=f'_{1}\circ f'^{-1}_{1}\circ f'_{2}(L'_{0})=f'_{2}(L'_{0})=B$. Hence, by induction if $\{A_{i}\}_{i=1}^{n}$ is a finite collection of closed smooth $G$-submanifolds of $X$ then there exists a smooth equivariant triangulation $f:K\rightarrow X$ and a finite collection of $G$-subcomplexes  $\{L_{i}\}_{i=1}^{n}$ of $K$ such that $L_{i}$ is a  $G$-triangulation of $A_{i}$, for $i=1,\dots,n$.

To realize $A$ as the fixed point set of some $h:X\rightarrow X$, it is necessary to remove every fixed point $x\in X$ of $f:X\rightarrow X$ outside of $A$. Hence, we need to extend the notion of {\it neighborhood by-passed} for a closed subset $A$ as in \cite[Definition 2.1]{s} in order to handle these undesired fixed points.. Thus, a $G$-invariant subset $A$ is said to be {\it $G$-neighborhood by-passed} if there exists an invariant open subset $U\subset X$ such that $A\subset U$ and $U$ can be by-passed in $X$.

We observe that if $(X,A)$ is a $G$-ENR pair then $A$ is an invariant neighborhood retract in $X$ and if $\{Y_{i}\}_{i=1}^{n}$ is a finite collection of closed smooth $G$-submanifolds such that $\dim(Y_{i})+1<\dim(X)$, then $A\cup Y$ ($Y=\bigcup_{i=1}^{n} Y_{i}$) remains a by-passed $G$-subset of $X$ provided $A$ is by-passed in $X$. Furthermore, a close inspection of the proof of Theorem 2.2 of \cite{s} indicates that the same argument works for the same result in the equivariant setting. That is, if $A$ is a by-passed locally contractible $G$-subset of $X$ then $A$ is $G$-neighborhood by-passed, for $X$ a compact smooth $G$-manifold with $\dim(X)\geq 3$. To see that, we note that if $K$ is the $G$-triangulation of $X$ then there is a by-passed neighborhood (may not be equivariant) $U$ of $A$ in $|K|$. We obtain the open $G$-subpolyhedron:

$$St(A,K)= \bigcup_{\stackrel{\overline{|t|}\cap A\neq \emptyset}{t\in K}} |t|,$$

\noindent such that $\overline{St(A,K)}$ is a subset of $U$ by taking a $G$-refinement $K'$ of $K$ if necessary, where $t$ is a simplex of $K$. Therefore, if $p:I\rightarrow X$ is a path with endpoints in $U-\overline{St(A,X)}$ and outside $\overline{St(A,X')}$ then using Corollary 3.3.11 of \cite{sp} we deform $p$ out of $\overline{St(A,X')}$.

Thus, if $\{Y_{i}\}_{i=1}^{n}$ is a finite collection of closed smooth $G$-submanifolds such that $\dim(Y_{i})+1<\dim(X)$ (thus each $Y_i$ has codimension at least $2$ in $X$ so that $Y_i$ can be by-passed in $X$), then $A\cup Y$ ($Y=\bigcup_{i=1}^{n} Y_{i}$) remains a by-passed $G$-subset of $X$ using a finite collection of $G$-subcomplex  $\{L_{i}\}_{i=1}^{n}$ of $K$ such that $L_{i}$ is a  $G$-triangulation of $Y_{i}$, for $i=1,\dots,n$.

The next lemma shows how the fixed points outside $A$ may be removed (see also \cite{fo}).

\begin{lema}\label{Proced} Let $\{Y_{i}\}_{i=1}^{n}$ be a finite collection of closed $G$-submanifolds of the  $G$-manifold $X$ such that $\dim(Y_{i})+1<\dim(X)$ and the action of $G$ outside $Y=\bigcup_{i=1}^{n} Y_{i}$ is free, where $G$ is a finite group. Let $f:X\rightarrow X$ be a $G$-selfmap, $A$ be a non-empty closed locally contractible and by-passed $G$-subset of $X$ such that $A\subset {\rm Fix}(f)$, there are no fixed points of $f$ in $Y-A$, and $f$ has a finite number of fixed points in $X-(A\cup Y)$. Let $x_{0}$ and $x_{1}$ be two fixed points of $f$ that are $G$-Nielsen equivalent from different orbits such that $x_{0}\in X-(A\cup Y)$ and $x_{1}\in X-(A\cup Y)$ or $x_{1}\in \partial(A)$, where $\partial(A)$ is the boundary of $A$ in $X$ and $q:I\rightarrow X$ a path with end points $q(0)=x_{0}$ and $q(1)=x_{1}$ such that $f\circ q$ is homotopic to $q$ relative to the endpoints. 

Then, $f$ is $G$-homotopic, relative to $(A\cup Y)$, to a $G$-selfmap $h: X\rightarrow X$ such that ${\rm Fix}(h)={\rm Fix}(f)-G\{x_{0}\}$.
\end{lema}

{\bf Proof of Lemma \ref{Proced}: } Since $A$ is locally contractible and can be by-passed in $X$, the discussion above shows that $A$ is $G$-neighborhood by-passed in $X$. Furthermore, $A\cup Y$ can be by-passed in $X$. Thus, the path $q$ is homotopic, relative to endpoints, to a path $q'(t)$ such that for $0\le t<1$, $q'(t)\in X-(A\cup Y)$ with $q'(0)=x_0, q'(1)=x_1$. Since $G$ acts freely on $X-Y$ and hence on $X-(A\cup Y)$, taking the $G$-translates of $q'$ yields $|G|$ paths from the orbit $G\{x_0\}$ to the orbit $G\{x_1\}$. Note that the segements $G\{q'([0,1))\}$ are disjoint while $\{G\{q'(1)\}\}$ consists of $[G:G_{x_1}]$ distinct endpoints. Here, the isotropy subgroup $G_{x_1}$ at $x_1$ is trivial if $x_1\in X-(A\cup Y)$. Now we coalesce these two fixed orbits in the same fashion as in \cite[Lemma 3.1]{w2}. (For slightly more general spaces in which normal arcs are used, see \cite[Theorem 2]{fo}.)

\hfill$\Box$

We will prove Theorem \ref{Teo;sod} before Theorem \ref{Cor,sod} and for the same reason we prove Theorem \ref{Teo;sod} by first establishing Lemma \ref{ida} and Lemma \ref{volta}.

\begin{teo}\label{Teo;sod} Let $G$ be a compact Lie group, $X$ be a compact smooth $G$-manifold and $A$ be a nonempty, closed, locally contractible $G$-subset of $X$ such that for each finite $WK$ we assume that $\dim(X^{K})\geq 3$, $\dim(X^{K})-\dim(X^{K}-X_{K})\geq2$ and $A^{K}$ is by-passed in $X^{K}$, for all $(K)\in {\rm Iso}(X)$. Then, given a $G$-map $f:X\rightarrow X$ there exists a $G$-map $h:X\rightarrow X$ $G$-homotopic to $f$ with ${\rm Fix}(h)=A$ if, and only if, the conditions $(C_G1)$ and $(C_G2)$, given in Theorem \ref{Cor,sod}, hold for $f$ relative to $A$.
\end{teo}

\begin{lema}\label{ida} Let $G$ be a compact Lie group, $X$ be a $G$-space $G$-ANR and $A$ be a nonempty closed $G$-subset of $X$. If $f:X\rightarrow X$ is a $G$-map $G$-homotopic to $h:X\rightarrow X$ such that ${\rm Fix}(h)=A$ then the conditions $(C_G1)$ and $(C_G2)$ given by Theorem \ref{Cor,sod} hold for $f$ relative to $A$.
\end{lema}

{\bf Proof of Lemma \ref{ida}:} Let $H:X\times I \rightarrow X$ be a $G$-homotopy which starts at $f$ and ends at $h$. Then $\overline{H}=H|_{(X\times\{0\})\cup(A\times I)}:(X\times\{0\})\cup(A\times I)\rightarrow X$ satisfies $(C_G1)$. If $F$ is a $WK$-essential fixed point class of $f^{K}$, then, there exists a path $p:I\rightarrow X^{K}$ such that $p(0)\in F$ and $p(1)\in J$, where $J\subset A^{K}$ is a $WK$-essential fixed point class of $h^{K}$, $H^{K}$-related to $F$ and $\{p(t)\}\sim\{\overline{H}^{K}(p(t),t)\}$. In fact, 
$$\{\overline{H}^{K}(p(t),t)\}\sim\underbrace{\{\overline{H}^{K}(p(t),0)\}}_{=\{f\circ p(t)\}}*\{\overline{H}^{K}(p(1),t)\}.$$
So, $(C_G2)$ is satisfied.

\hfill$\Box$

Lemma \ref{ida} shows that the conditions $(C_G1)$ and $(C_G2)$ are necessary for $A={\rm Fix}(h)$. The example below shows that these two conditions are independent of each other.

\begin{exem}\label{exe2} Let $G=\mathbb{Z}_{2}$, $X=\mathbb{S}^{2}$ and the action given by $\xi(x,y,z)\mapsto(-x,-y,z)$. Then, there is no $\mathbb{Z}_{2}$-homotopy $H$ from the identity $Id$ to $h$ such that ${\rm Fix}(h)=\{(x,y,0)\in \mathbb{S}^{2}\}$. Note that $(C_G1)$ occurs, because the map is the identity, but $(C_G2)$ does not. On the other hand, let $G=\mathbb{Z}_{2}$, $X=\mathbb{S}^{3}$ and the action given by $\xi(x,y,z,w)\mapsto(x,y,z,-w)$. Then, there is no $\mathbb{Z}_{2}$-homotopy $H$ from the antipodal map $-Id$ to $h$ such that ${\rm Fix}(h)=\{(x,y,z,0)\in \mathbb{S}^{3}\}$. This time $(C_G2)$ holds because the map is fixed point free but $(C_G1)$ does not hold.
\end{exem}

\begin{lema}\label{volta} Let $G$ be a compact Lie group, $X$ be a compact smooth $G$-manifold and $A$ be a nonempty, closed, locally contractible $G$-subset of $X$ such that for each finite $WK$ we assume that $\dim(X^{K})\geq 3$, $\dim(X^{K})-\dim(X^{K}-X_{K})\geq2$ and $A^{K}$ is by-passed in $X^{K}$, for all $(K)\in {\rm Iso}(X)$. If the conditions $(C_G1)$ and $(C_G2)$, given in Theorem \ref{Cor,sod}, hold for a $G$-map $f:X\rightarrow X$ relative to $A$, then there exists a $G$-map $h:X\rightarrow X$, $G$-homotopic to $f$ with ${\rm Fix}(h)=A$.
\end{lema}

{\bf Proof of Lemma \ref{volta}:} This proof follows the steps of the proof of Theorem 3.2 of \cite{sc1}. Consider a $G$-map $\overline{H} : (X\times\{0\})\cup(A \times I) \rightarrow X$ given by $(C_G1)$. It is possible to extend $\overline{H}$ to a $G$-homotopy $_{1}\overline{H}_{1} : \Big( X\times \{0\} \Big) \cup \Big( (A\cup X_{1}) \times I \Big) \rightarrow X$. As commented above, there is a closed $G$-invariant neighborhood $V$ of $A_{1}$ inside $X_{1}$ and $V$ retracts onto $A_{1}$ equivariantly. Note that $WK_{1}$ acts freely on $X^{K_{1}}_{1}=X_{K_{1}}$ and $_{1}h^{K_{1}}_{1}$ is a $WK_{1}$-map. Hence, if $WK_{1}$ has positive dimension we apply Lemma 3.3 of \cite{da} and Lemma 2.1 of \cite{fw} to extend $_{1}\overline{H}_{1}$ to a $G$-homotopy $\overline{H}_{1}:(X\times\{0\})\cup((A\cup X_{1})\times I)\rightarrow X$, relative to $V$. Moreover, $h_{1}$ has no fixed points in $X_{1}-A_{1}$ and ${\rm Fix}(h_{1})=A$, where $h_{1}= \overline{H}_{1}(\bullet,1):A\cup X_{1} \rightarrow X$.

On the other hand, if $WK_{1}$ is a finite group then $X^{K_{1}}$ is a $WK_{1}$-polyhedron such that $A_{1}^{K_{1}}$ is a $WK_{1}$-subpolyhedron and $St(A_{1},X^{K_{1}})$ is neighborhood by-passed in $X^{K_{1}}$. We apply Lemma 3.1 of \cite{da} and Lemma \ref{Proced} to obtain a $WK_{1}$-homotopy $H :(A_{1}\cup X_{1})^{K_{1}} \times I  \rightarrow X^{K_{1}}$ which can be extended by Lemma 2.1 of \cite{fw} to a $G$-homotopy $\overline{H}_{1}:(X\times\{0\})\cup((A\cup X_{1})\times I)\rightarrow X$, relative to $V$, such that $h_{1}$ has no fixed points in $X_{1}-A_{1}$ and ${\rm Fix}(h_{1})=A$, where $h_{1}= \overline{H}_{1}(\bullet,1):A\cup X_{1} \rightarrow X$.

By induction, we may assume that we have a $G$-map $\overline{H}_{i-1} : (X\times \{0\})\cup ((A\cup X_{i-1}) \times I) \rightarrow X$ such that ${\rm Fix}(h_{i-1})=A$, where $h_{i-1}= \overline{H}_{i-1}(\bullet,1):A\cup X_{i-1} \rightarrow X$ and the proof follows the steps we did for $WK_{1}$.

\hfill$\Box$

Now Theorem \ref{Teo;sod} follows easily from Lemma \ref{ida} and Lemma \ref{volta}.

{\bf Proof of Theorem \ref{Cor,sod}:} First of all, by Theorem \ref{Teo;sod}, there is a $G$-map $h_{1}: X\rightarrow X$ $G$- homotopic to $f$ such that ${\rm Fix}(h_{1})=A$. We may apply Proposition 2.5 of \cite{da} and Theorem 4.3 of \cite{w1} to conclude that $h_{1}$ is $G$-homotopic to $h_{2}$ such that $h_{2}|_{X^{K}}$ has a finite number of fixed points, all of which inside $St(A^{K})$ and lying in the interior of a maximal simplex of $X^{K}$ and $h_{2}$ is a $G$-proximity map in $St(A)$ (for some $G$-triangulation of $X$).

Since $\Phi$ has nonempty intersection with every component of $A$ we can pull the fixed points of $h_{2}$ to $\Phi$. Let $\alpha$ be the $G$-map of Lemma VIII.C.1 of \cite{br} and $\overline{d}$ the equivariant bounded distance in $X$ then we define 
$$\overline{H}_{3}:(X\times\{0\})\cup(St(A)\times I) \rightarrow X $$
given by:

\begin{displaymath}
\begin{array}{ccc}
   (x,t) & \mapsto & \left\{ \begin{array}{ll}
                                        \alpha(x,h_{2}(x),1 - (1-\overline{d}(x,\Phi))t)& if\ (x,t)\in St(A)\times I; \\
                                         h_{2}(x)                                    & if\ t=0 .
                                         \end{array}\right.
\end{array}
\end{displaymath}

Then, we extend $\overline{H}_{3}$, relative to $\partial(St(A))$, to a $G$-map $\overline{H}_{4}:X\times I\rightarrow X$. By Lemma 3.1 of \cite{da}, we eliminate the fixed points of $\overline{H}_{4}(\bullet,1)$ inside $X-St(A)\times\{1\}$. This finite set of fixed points can be removed because these fixed points lie in some non essential fixed point classes of $\overline{H}_{4}(\bullet,1)$ since $h_2|_{X-Int(A)}$ is fixed point free. Thus, the resulting $G$-map is a $G$-homotopy $H: X\times I \rightarrow X$  connecting $f$ to a $G$-map $h$ such that ${\rm Fix}(h)=\Phi$.

%


%

\section[Proof of Theorem 1.6]{Proof of Theorem \ref{Teo;bs} - An equivariant analogue of a theorem of Brown-Soderlund}\label{sec3}

Throughout this last section, $G$ will denote a finite group.
Given a $G$-fiber-preserving map $f:X\rightarrow X$ of the total space $X$ of a $G$-fiber bundle $\mathfrak{F}=(X,p,B,Y)$, it is known that the fixed point set of $f$ is related with the fixed point set of the induced map $\overline{f}:B\rightarrow B$. However, there are equivariant homotopies that are not fiber-preserving as in the example below:

\begin{exem}\label{Exe;s3}
Let $G=\mathbb{Z}_{2}$ and $X=\mathbb{S}^{2}\times \mathbb{S}^{1}$ and the action is given by $\xi((a,b,c), \cos x+i\sin x)\mapsto((a,b,c), \cos x-i\sin x)$. The $G$-map $f$, defined on $X$ by setting $f((a,b,c), \cos x+i\sin x)=((-a,-b,-c), \cos x+i\sin x)$, is the start of the following equivariant homotopy:

\begin{displaymath}
\begin{array}{c}
H\Big( \big( \cos\theta\sin\psi {\bf,} \sin\theta\sin\psi ,  \cos\psi \big), \cos x+i\sin x,t \Big)= \\
\Big( \big(- \cos(\theta+t|\sin x |\pi)\sin\psi, -\sin(\theta+t|\sin x |\pi)\sin\psi ,- \cos\psi \big),
\cos x+i\sin x \Big) .
\end{array}
\end{displaymath}

Then, $A=\{(a,b,0)\in \mathbb{S}^{2}\}\times\{-i,i\}=\mathbb{S}^{1}\times\{i,-i\}$ is the fixed point set of $h\in[f]_G$ where $h=H(\bullet,1)$. Let $p=\pi_{1}:\mathbb{S}^{2}\times \mathbb{S}^{1}\rightarrow \mathbb{S}^{2}$ be the projection, then $(\mathbb{S}^{2}\times \mathbb{S}^{1},\pi_{1},\mathbb{S}^{2})$ is a $\mathbb{Z}_{2}$-fiber bundle, $f$ is a fiber-preserving map and the induced map $\overline{f}=a:\mathbb{S}^{2}\rightarrow \mathbb{S}^{2}$ is the antipodal map. However, $p((x,y,z), 1)=(x,y,z)= p((x,y,z), i)$ and $p\circ h((x,y,z),1) =(-x,-y,-z)$ is different from $p\circ h((x,y,z),i) = (x,y,-z)$. So, $h$ is not a fiber preserving map and $H$ is not a fiber-preserving homotopy. In fact, $A$ {\it cannot} be realized as the fixed point set of {\it any} map equivariantly fiberwise homotopic to $f$. To see that, we note that $X^G=\mathbb S^2 \times \{\pm 1\}=\mathbb S^2_1 \sqcup \mathbb S^2_{-1}$, where $(w, \pm 1)\in \mathbb S^2_{\pm 1}$, consists of two disjoint $2$-spheres $\mathbb S^2$. If $F_t$ is a $\mathbb Z_2$ fiber-preserving homotopy such that $F_0=f$ and ${\rm Fix}(F_1)=A$, then $F^G_t$ is a homotopy on $X^G$. Now, $f^G=F^G_0$ maps $\mathbb S^2_1$ to $\mathbb S^2_1$ and $\mathbb S^2_{-1}$ to $\mathbb S^2_{-1}$. On the other hand, $F_1$ is fiber-preserving and $A$ is the fixed point of $F_1$, it follows that the induced map $\overline {F_1}$ fixes the circle $\{(a,b,0)\in \mathbb S^2\}$ pointwise. This implies that  $F_1$ maps the (non-fixed) point $((a,b,0),1)$ to the point $((a,b,0),-1)$ so that $F_1$ maps the equator of $\mathbb S^2_{1}$ to that of $\mathbb S^2_{-1}$, and vice versa. Thus $F^G_1$ maps $X^G$ to itself by interchangeing the two disjoint spheres $\mathbb S^2_{\pm 1}$. The images of $X^G$ under $F^G_0$ and $F^G_1$ contradict the continuity of $F^G_t$. Hence such an equivariant fiber-preserving homotopy $F_t$ cannot exist.    
\end{exem}

The example above indicates the importance of modifying the conditions $(C_G1)$ and $(C_G2)$ and replacing them by $(C_G1)_{\mathfrak{F}}$ and $(C_G2)_{\mathfrak{F}}$  for the fiber-preserving map setting.

\begin{lema}\label{lema3:GF} Let $f:X\rightarrow X$ be a $G$-fiber preserving map in the total space of the $G$-fiber bundle $\mathfrak{F}=(X,p,B,Y)$, where $X$, $B$ and $Y$ are $G$-spaces $ANR$. Suppose that there is a $G$-fiber preserving homotopy connecting a $G$-fiber preserving map $h:X\rightarrow X$ to $f$ such that ${\rm Fix}(h)=A$ for a nonempty and closed $G$-subset $A$ of $X$. Then the conditions $(C_G1)_{\mathfrak{F}}$ and $(C_G2)_{\mathfrak{F}}$ given in Theorem \ref{Teo;bs} hold for $f$ and $A$.
\end{lema}

The proof of Lemma \ref{lema3:GF} follows the steps of Lemma \ref{ida}. Since $\mathfrak{F}=(X,p,B,Y)$ is a $G$-fiber bundle where $X$, $B$ and $Y$ are compact smooth $G$-manifolds, we observe that $(X,p,B)$ is a $G$-fibration and there is a $G$-lift map $\Lambda:\Omega_{p}\rightarrow E^{I}$ such that $\Lambda(e,\alpha)(0)=e$, $p\circ\Lambda(e,\alpha)(t)=\alpha(t)$ and $\Lambda(e,p(e))(t)=e$, for all $t\in I$, where $E^{I}=\{\alpha:I\rightarrow E; \ \alpha \ is \ a \ path \}$ and $\Omega_{p}=\{(e,\alpha)\in X\times B^{I}; \ p(e)=\alpha(0) \}$. 

\begin{rem} We should point out that Lemma \ref{lema3:GF} holds for any compact Lie group $G$ if we modify condition $(C_G2)_{\mathfrak{F}}$ by only considering those $(K)$'s with $|WK|<\infty$.
\end{rem}

The next proposition is an equivariant analogue of Theorem 2.1 of \cite{af}.

\begin{prop}\label{prop2:GF} Let $\overline{H}:(X\times\{0\})\cup (A\times I)\rightarrow E$ be a $G$-map in the $G$-fibration $\mathfrak{F}=(E,p,B)$, where $E$ is a $G$-ANR,  $A$ is a closed $G$-subset of $X$, $(X,A)$ is a $G$-metric pair and $p\circ \overline{H}(x,0)=p\circ {\overline{H}}(x,t)$ for all $(x,t)\in A\times I$. Then $\overline{H}$ can be extended to a $G$-homotopy $H:X\times I\rightarrow E$ such that $p\circ H(x,0)=p\circ H(x,t)$ for all $(x,t)\in X\times I$.
\end{prop}

{\bf Proof of Proposition \ref{prop2:GF}:} Let $H':X\times I\rightarrow E$ a $G$-extension of $\overline{H}$. Then $H'$ is given by:

\begin{displaymath}
	\begin{array}{ccclccl}
	H': & X & \rightarrow & E^{I}            &   &             &    \\
	    & x & \mapsto     & H'(x,\bullet):   & I & \rightarrow & E  \\
      &   &             &                  & t & \mapsto     & H'(x,t).  	
  \end{array}
\end{displaymath}

Then define $H(x,t)=\Lambda(H'(x,t),p(H'(x,\bullet))_{t})(1)$, where $p(H'(x,\bullet))_{t}(s)=p(H'(x,(1-s)t))$ and $\Lambda$ is a $G$-lift map.

\hfill$\Box$

\begin{lema}\label{lema1:GF} 
Let $\mathfrak{F}=(X,p,B,Y)$ be a $G$-fiber bundle where $X,B$ and $Y$ are compact and smooth $G$-manifolds, $\dim(B^{K})\geq 3$, $\dim(B^{K})-\dim(B^{K}-B_{K})\geq2$, for all $(K)\in {\rm Iso}(B)$, $A$ be a nonempty, closed, locally contractible $G$-subset of $X$ such that $p(A)$ be a closed $G$-subset of $B$ and $p^{K}(A^{K})$ is by-passed in $B^{K}$, for all $(K)\in {\rm Iso}(B)$, and $f:X\rightarrow X$ a $G$- fiber preserving map such that conditions $(C_G1)_{\mathfrak{F}}$ and $(C_G2)_{\mathfrak{F}}$ given in Theorem \ref{Teo;bs} hold for $f$ and $A$.

Then there exists a $G$-fiber-preserving map $h$, $G$-fiberwise homotopic to $f$ with $A\subset {\rm Fix}(h)\subset p^{-1}(p(A))$ and ${\rm Fix}(\overline{h})\cap (B-p(A))$ is a finite set.
\end{lema}

{\bf Proof of Lemma \ref{lema1:GF}:} $p(A)$ is a closed $G$-subset of $B$ then the $G$-fiber-preserving map $H_{A}:A\times I\rightarrow X$ given by $(C_G1)_{\mathfrak{F}}$ induces a $G$-map $\overline{H}_{A}:p(A)\times I\rightarrow B$ such that $\overline{H}_{A}(\bullet,0)=\overline{f}$ and $\overline{H}_{A}(\bullet,1)=i_{p(A)}:p(A)\hookrightarrow B$ the inclusion map.

Observe that we have almost the same conditions that we had in Theorem \ref{Cor,sod} except for $(C_G2)$. In this situation, suppose we have a $G$-map $\overline{H}_{i-1,A}:(p(A)\cup B_{i-1})\times I\rightarrow B$. As commented in Lemma \ref{volta}, it is possible to extend $\overline{H}_{i-1,A}$ to a $G$-map $\overline{H}_{i,1}:(B_{i}\cup p(A))\times I\rightarrow B$ relative to $p(A)\cup B_{i-1}$.

Since $WK_{i}$ is a finite group, $B_{i}^{K_{i}}$ is a $WK_{i}$- polyhedron such that $B_{i-1}^{K_{i}}$ is a $WK_{i}$-subpolyhedron of $B_{i}^{K_{i}}$ and $St(p(A_{i}^{K_{i}}))$ is neighborhood by-passed in $B_{i}^{K_{i}}$. Let $V$ be a $G$-invariant neighborhood retract of $St(p(A_{i}))\cup B_{i-1}$. It follows from Lemma 3.1 of \cite{da} and Lemma \ref{Proced} that there exists of a $WK_{i}$-homotopy $\overline{H}_{i}:B_{i}^{K_{i}}\times I\rightarrow B_{i}^{K_{i}}$ from $\overline{H}_{i,1}^{K_{i}}(\bullet,1)$ to $\overline{h}=\overline{H}_{i}(\bullet,1)$ such that:

\begin{enumerate}
\item $p(A)^{K_{i}}_{i}\subset {\rm Fix}(\overline{h})$;

\item $\overline{h}$ has a finite number of fixed points in $B^{K_{i}}_{i}-V^{K_{i}}$;

\item given a $WK_{i}$-fixed point class $F$ of $\overline{h}$ such that $F\cap p(A)^{K_{i}}_{i}=\emptyset$ then $F=WK_{i}\{x\}$, where $x\in B^{K_{i}}_{i}-V^{K_{i}}$ and $F$ is an essential $WK_{i}$-fixed point class of $\overline{h}$.
\end{enumerate}

Then, the $G$-map given by:
 \begin{displaymath}
	\begin{array}{rcl}
	\overline{h}_{t}(x) & = & \left\{\begin{array}{ll}
     	          g\overline{H}_{i}(g^{-1}x,t), & for\ x\in X-A, \ where \ G_{x}=g WK_{i} g^{-1};\\
                \overline{h}(x),                & for\ x\in V.
                        \end{array}\right.
	\end{array}
\end{displaymath}
extends a $WK_{i}$-homotopy to a $G$-homotopy $\overline{H}_{i}:(B_{i}\cup p(A))\times I\rightarrow B$ relative to $V$ and such that  
$${\rm Fix}(\overline{H}_{i}(\bullet,1))= p(A)\cup \big (\bigcup_{j\in \ T, \ j\leq \ i}(G\{b_{j,1}\}\cup\dots\cup G\{b_{j,m_{j}}\}) \ \big )$$
and $WK_{i}\{b_{i,l}\}$ is a essential $WK_{i}$-fixed point class of $\overline{H}_{i}^{K_{i}}(\bullet,1)$, for $1\leq l\leq m_{i}$.

Observe that if $p^{K}(F)=WK\{b_{i,l}\}$ for an essential $WK$-fixed point class $F$ of $f^{K}$ where $(K)\in {\rm Iso}(X)$, then we have a path $\overline{\alpha}$ such that:
$$\{\overline{\alpha}\}\thicksim \{\overline{f}^{K}\circ \overline{\alpha}\}*\{\overline{H}^{K}_{A}(\overline{\alpha}(1),t)\}\thicksim \{\overline{H}^{K}(\overline{\alpha}(t),t)\}.$$

\noindent Hence, $\overline{\alpha}(1)=gb_{i,l}$, for some $g\in WK$ and $\overline{\alpha}(1)\in p^{K}(A)$. However, this cannot occur because $b_{i,l}\notin p^{K}(A)$ and $p^{K}(A)$ is $WK$-invariant. By induction we extend the $G$-map $\overline{H}_{A}:p(A)\times I\rightarrow B$ to a $G$-homotopy $\overline{H}:B\times I\rightarrow B$ with the properties above.

Note that $H':X\times I\rightarrow B$ defined by $H'(x,t)=\overline{H}(p(x),t)$ is such that 
$$H'(x,0)=\overline{H}(p(x),0)=\overline{f}\circ p(x)=p\circ f(x).$$

\noindent Therefore, the lift of $H'$ is a fiber-preserving $G$-homotopy $H_{1}:X\times I\rightarrow X$ such that $f(x)=H_{1}(x,0)$ and $h_{1}(x)=H_{1}(x,1)$. Thus,
$${\rm Fix}(h_{1})\subset p^{-1}({\rm Fix}(\overline{h}))=p^{-1}(p(A)\cup G\{b_{1}\}\cup\dots\cup G\{b_{l}\}). $$

For each $G$-orbit $G\{b_{j}\}$ take the restriction $h_{1,b_{j}}$ of $h_{1}$ for $Gp^{-1}(b_{j})=p^{-1}(G\{b_{j}\})$, so $h_{1,b_{j}}:Gp^{-1}(b_{j})\rightarrow Gp^{-1}(b_{j})$ has no essential fixed point classes. In fact, suppose that $h_{1,b_{j}}^{K}$ has an essential $WK$-fixed point class $F$. Then, given $x\in F$ we have $WK\{x\}$ lying inside an essential $WK$-fixed point class of $h_{1}^{K}$. Thus, there exists a  $WK$-fixed point class $Q$ of $h_{1}^{K}$ which contains $WK\{x\}$. But, $h_{1}^{K}$ is fiber-preserving $WK$-homotopic to $f^{K}$, so, there exists an essential $WK$-fixed point class $D$ of $f^{K}$ $H_{1}^{K}$-related to $Q$. Note that $p^{K}(D)$ cannot be $\overline{H}^{K}$-related to $WK\{b_{j}\}$. Consequently, $h_{1,b_{j}}$ is fiber-preserving $G$-homotopic to $h_{2,b_{j}}: Gp^{-1}(b_{j})\rightarrow Gp^{-1}(b_{j})$ fixed point free.

Consider the $G$-map 
$$\widetilde{H}_{2}:(X\times \{0\}) \cup(p^{-1}({\rm Fix}(\overline{h}))\times I)\rightarrow X$$

\noindent defined by:

\begin{displaymath}
	\begin{array}{lcr}
	 \widetilde{H}_{2}(x,t) & = & \left\{ \begin{array}{ll}
                                        h_{1}(x)         & if\ t=0\ or\ if\ x\in p^{-1}(p(A)) ;\\
                                        H_{2,b_{j}}(x,t) & if\ x\in p^{-1}(G\{b_{j}\}).
                                         \end{array}\right.
\end{array} 
\end{displaymath}

\noindent With Proposition \ref{prop2:GF} we extend $\widetilde{H}_{2}$ to a fiber-preserving $G$-homotopy $H_{2}:X\times I\rightarrow X$ and $h_{2}=H_{2}(\bullet,1)$ is such that $\overline{H}_{2}(\bullet,1)=\overline{h}$. By $(C_G1)_{\mathfrak{F}}$, $h|_{2}$ is fiber-preserving $G$-homotopic to $i_{A}$. Let $\widetilde{H}_{A}$ such that $h_{2}|_{A}=\widetilde{H}_{A}(\bullet,0)$ and $i_{A}=\widetilde{H}_{A}(\bullet,1)$. Define $\widetilde{H}: (X\times \{0\})\cup((A \cup p^{-1}(G\{b_{1},\dots,b_{r}\}))\times I)  \rightarrow  X$ given by:

\begin{displaymath}
\begin{array}{ccl}
  \widetilde{H} (x,t) &  = & 
     \left\{\begin{array}{ll}
     h_{2}(x),               &if\ t=0;\\
     h_{2,b_{j}}(x),         &if\ x\in Gp^{-1}(b_{j});\\
     \widetilde{H}_{A}(x,t), &if\ x\in A.
     \end{array}\right.
\end{array}
\end{displaymath}

\noindent Applying Proposition \ref{prop2:GF} again we extend $\widetilde{H}$ to a fiber-preserving $G$- homotopy $H:X\times I\rightarrow X$ such that $A\subset {\rm Fix}(h)\subset p^{-1}(p(A))$ and ${\rm Fix}(\overline{h})\cap (B-p(A))$ is a finite set.

\hfill$\Box$

\begin{lema}\label{lema2:GF}
Let $(\mathfrak{F},\mathfrak{F}_{0})=((X,A),p,B,(Y,Y_{0}))$ be a $G$-fiber bundle pair, where $X$, $B$ and $Y$ are compact and smooth $G$-manifolds, $B$ retracts equivariantly to a point $b_{0}\in B$ and $\dim(Y^{K})\geq 3$ and $\dim(Y^{K})-\dim(Y^{K}-Y_{K})\geq2$, for all $(K)\in {\rm Iso}(Y)$. Let $Y_{0}$ be a closed and locally contractible $G$-subset of $Y$ such that $Y^{K}_{0}$ is by-passed in $Y^{K}$, for all $(K)\in {\rm Iso}(Y)$, $A$ be a nonempty, closed, locally contractible $G$-subset of $X$ and $f:X\rightarrow X$ be a $G$-map such that $p\circ f=p$, $A\subset {\rm Fix}(f)$, $A^{K}$ intersects every essential $WK$-fixed point class of $f_{b_{0}}^{K}:WK(p^{K})^{-1}(\{b_{0}\})\rightarrow WK(p^{K})^{-1}(\{b_{0}\})$, for all $(K)\in {\rm Iso}(X)$.

Then for every closed $G$-invariant subset $Z$ of $A$ that intersects every component of $A$ and $(A,Z)$ is $G$-fiber bundle pair of $\mathfrak{F}_{0}$ there exists a fiber-preserving $G$-map $h$, $G$-fiberwise homotopic to $f$ with ${\rm Fix}(h)=Z$. 
\end{lema}

{\bf Proof of Lemma \ref{lema2:GF}:} $(X,p,B)$ is $G$-equivalent to a trivial $G$-fibration  $(B\times Y,\pi,B)$, where $\pi$ is a projection in $B$. So, there exists a $G$- homeomorphism $\Phi:B\times Y\rightarrow X$ such that $\Phi(B\times Y_{0})=A$ and $p\circ\Phi=\pi$. Define $f^{*}=\Phi^{-1}\circ f\circ \Phi:B\times Y\rightarrow B\times Y$ and note that:

$$\pi\circ f^{*}= \underbrace{(p\circ\Phi)\circ (\Phi^{-1}}_{=p}\circ f\circ \Phi) =\underbrace{p\circ f}_{=p}\circ \Phi=p\circ\Phi=\pi.$$

\noindent Therefore, $f^{*}(b,y)=(b,f^{*}_{b}(y))$ and $gf^{*}(b,y)=(gb,f^{*}_{gb}(gy))$, for all $g\in G$.

$B$ retracts equivariantly to $b_{0}$, so, there exists a $G$-homotopy $D:B\times I\rightarrow B$ such that for each $b\in B$ we have $D(b,0)=b$ and $D(b,1)=b_{0}$. Then, define $U^{*}:B\times Y\times I\rightarrow B\times Y$ given by:

\begin{displaymath}
U^{*}(b,y,t)=\left\{\begin{array}{lc}
(b,f^{*}_{D(b,2t)}(y)),       & if\ 0\leq t\leq \frac{1}{2} \\
(b,f^{*}_{D(b_{0},2-2t)}(y)), & if\ \frac{1}{2}\leq t\leq 1.\\
\end{array}\right.
\end{displaymath}

\noindent Note that $U^{*}(b,y,0)=(b,f^{*}_{b}(y))=f^{*}(b,y)$ and $U^{*}(b,y,1)=(b,f^{*}_{b_{0}}(y))$. Then, $f^{*}$ is $G$-homotopic to $id\times f^{*}_{b_{0}}$. Since $\Phi(B\times Y_{0})= A\subset {\rm Fix}(f)$ we have, for each $(b,y)\in B\times Y_{0}$:

$$f^{*}(b,y)=\Phi^{-1}\circ \underbrace{f\circ\underbrace{\Phi(b,y)}_{\in A}}_{=\Phi(b,y)}=\Phi^{-1}\circ\Phi(b,y)=(b,y).$$

\noindent Then, $B\times Y_{0}\subset {\rm Fix}(f^{*})$ and $Y_{0}\subset {\rm Fix}(f^{*}_{b})$ because $(b,f^{*}_{b}(y))=f^{*}(b,y)=(b,y)$. 

By hypothesis, $A$ intersects each essential $WK$-fixed point class of $f^{K}_{b_{0}}:WK(p^{K})^{-1}(b_{0}) \rightarrow WK(p^{K})^{-1}(b_{0})$. So, $A\cap WK(p^{K})^{-1}(b_{0})$ intersects each essential $WK$-fixed point class of $f_{b_{0}}^{K}$. So, $Y_{0}$ intersects each essential $WK$-fixed point class of $(f^{*}_{b_{0}})^{K}$ because:

$$(\Phi^{-1})^{K}(A^{K}\cap WK(p^{K})^{-1}(b_{0}))=WK\{b_{0}\}\times Y_{0}^{K}.$$

The $G$-fiber bundle pair $\Big((A,Z), p, B, (Y_{0},\Omega)\Big)$ is such that $\Omega$ intersects every component of $Y_{0}$ because $Z$ intersects every component of $A$. Therefore,  $(C_G1)_{\mathfrak{F}}$ and $(C_G2)_{\mathfrak{F}}$ hold for $Y_{0}$ and $f^{*}_{b_{0}}$. By Theorem \ref{Teo;sod} there exists a homotopy $V^{*}:Y\times I\rightarrow Y$ such that $f^{*}_{b_{0}}=V^{*}(\bullet,0)$, $V^{*}(\bullet,1)=g^{*}_{b_{0}}:Y\rightarrow Y$ and ${\rm Fix}(g^{*}_{b_{0}})=\Omega$. Define a fiber-preserving $G$-homotopy $H^{*}:B\times Y\times I\rightarrow B\times Y$ given by:

\begin{displaymath}
H^{*}(b,y,t)=\left\{\begin{array}{lc}
U^{*}(b,y,2t),     & if\ 0\leq t\leq \frac{1}{2} \\
(b,V^{*}(y,2t-1)), & if\ \frac{1}{2}\leq t\leq 1.\\
\end{array}\right.
\end{displaymath}

\noindent Therefore, $f^{*}$ is $G$-homotopic to $id\times g^{*}_{b_{0}}$ and ${\rm Fix}(id\times g^{*}_{b_{0}})=B\times\Omega$. Then, $H:X\times I\rightarrow X$ given by $H(e,t)=\Phi\circ H^{*}(\Phi^{-1}(e),t)$ is a $G$-homotopy such that $H(x,0)=f(x)$ and ${\rm Fix}(H(\bullet,1))=\Phi(B\times\Omega)=Z$.

\hfill$\Box$ 

\ %

{\bf Proof of Theorem \ref{Teo;bs}:} With Lemma \ref{lema1:GF} we assume that:

\begin{enumerate}
	\item $A\subset {\rm Fix}(f) \subset p^{-1}(p(A))$;
	
	\item $F={\rm Fix}(\overline{f})\cap (B-p(A))$ is a finite set.
\end{enumerate}

Let $f_{j}=f|_{p^{-1}(p(A)_{j})}:p^{-1}(p(A)_{j})\rightarrow p^{-1}(p(A)_{j}) $ a restriction of $f$, so $p\circ f_{j}=p$. Using $X=p^{-1}(p(A)_{j}), \ A=A\cap p^{-1}(p(A)_{j}), \ B=p(A)_{j}, \ b_{0}=b_{j}, \ Y_{0}=Y_{j}$ and $f=f_{j}$ the hypotheses of Lemma \ref{lema2:GF} are satisfied and there exists a fiber-preserving $G$-homotopy $H_{j}:p^{-1}(p(A)_{j})\times I\rightarrow p^{-1}(p(A)_{j})$ from $f_{j}$ to $h_{j}=H_{j}(\bullet,1)$ such that ${\rm Fix}(h_{j})=Z_{j}$.

Define $ \widetilde{H}_{2}:  (X\times \{0\})\cup(p^{-1}(F\cup p(A))\times I)  \rightarrow  X$ by:
\begin{displaymath}
\begin{array}{ccl}
 \widetilde{H}_{2} (x,t)  & = & 
   \left\{\begin{array}{ll}
           f(x),       & if\ t=0\ or\ p(x)\in F \\
           H_{j}(x,t), & if\ p(x)\in p(A)_{j}.
   \end{array}\right.
\end{array}
\end{displaymath}

\noindent With Proposition \ref{prop2:GF} there is a fiber-preserving $G$-homotopy $H:X\times I\rightarrow X$ such that $p(H(x,t))=\overline{f}\circ p(x)$. Therefore, $h=H(\bullet,1):X\rightarrow X$ is such that ${\rm Fix}(h)=Z$ and $h$ is fiber-preserving $G$-homotopic to $f$.

\hfill$\Box$

\begin{cor}\label{cor2:GF}  Let $\mathfrak{F}=(X,p,B,Y)$ be a $G$-fiber bundle where $X$, $B$ and $Y$ are compact and smooth $G$-manifolds, $\dim(B^{K})\geq 3$, $\dim(B^{K})-\dim(B^{K}-B_{K})\geq2$, for all $(K)\in {\rm Iso}(B)$, $\dim(Y^{K})\geq 3$, $\dim(Y^{K})-\dim(Y^{K}-Y_{K})\geq2$, for all $(K)\in {\rm Iso}(Y)$.

Let $A$ be a nonempty, closed, locally contractible $G$-subset of $X$ such that $(X,A)$ is $G$-fiber bundle pair with respect to the fiber bundle $\mathfrak{F}$, $p(A)$ be a closed $G$-subset of $B$ such that each component $p(A)_{j}$ of $p(A)$ is equivariantly contractible and $p^{K}(A^{K})$ is by-passed in $B^{K}$, for all $(K)\in {\rm Iso}(B)$. Let $Y_{j}$ be a subbundle fiber of $A$ such that $Y_{j}$ is a closed and locally contractible $G$-subset of $Y$ and $Y^{K}_{j}$ is by-passed in $Y^{K}$, for all $(K)\in {\rm Iso}(Y)$, and $f:X\rightarrow X$ be a $G$-fiber-preserving map such that $A^{K}$ intersects every essential $WK$-fixed point class of $f^{K}_{b_{j}}:WK(p^{K})^{-1}(\{b_{j}\})\rightarrow WK(p^{K})^{-1}(\{b_{j}\})$ for at least one $b_{j}$ in each component $p^{K}(A^{K})_{j}$, for all $(K)\in {\rm Iso}(X)$.

Then there exists a $G$-fiber-preserving map $h$, $G$-fiberwise homotopic to $f$ with ${\rm Fix}(h)=A$ if, and only if, the following conditions holds for $f$ and $A$:

\begin{itemize}
\item[$(C_G1)_{\mathfrak{F}}$]  there exists a $G$-fiber-homotopy $H_{A}:A\times I\rightarrow X$ from $f|_{A}$ to the inclusion $i:A\hookrightarrow X$;
\item[$(C_G2)_{\mathfrak{F}}$]  for every $WK$-essential fixed point class $F$ of $f^{K}:X^{K}\rightarrow X^{K}$ there exists a path $\alpha: I\rightarrow X^{K}$ with $\alpha(0)\in F$, $\alpha(1)\in A^{K}$, and $\{\alpha(t)\}\sim \{f^{K}\circ\alpha(t)\}\ast\{H_{A}^{K}(\alpha(1),t)\}.$
\end{itemize}
\end{cor}

{\bf Proof of Corollary \ref{cor2:GF}:} If the conditions hold then we apply Theorem \ref{Teo;bs} for $Z=A$. If there exists $h$ then by Lemma \ref{lema3:GF} $(C_G1)_{\mathfrak{F}}$ and $(C_G2)_{\mathfrak{F}}$ hold.

\hfill$\Box$

\end{document}